\documentclass[10pt,A4]{article}
\usepackage[latin5]{inputenc}
\usepackage{amsmath,amssymb}
\usepackage{amsthm}
\usepackage{indentfirst}

\newtheorem{theorem}{Theorem}

\newtheorem{corollary}{Corollary}

\newtheorem{proposition}{Proposition}
\newtheorem{remark}{Remark}
\numberwithin{definition}{section} \numberwithin{theorem}{section}
\numberwithin{lemma}{section}\numberwithin{corollary}{section}
\numberwithin{equation}{section} \numberwithin{example}{section}
\numberwithin{proposition}{section} \numberwithin{remark}{section}
\begin{document}

\begin{center}

{\bf \Large Approximations of the Set of  Trajectories \\ and Integral Funnel of the Nonlinear Control Systems \\ \vspace{2mm} with Limited Control Resources}

\vspace{3mm}

Nesir Huseyin$^1$, Anar Huseyin$^2$, Khalik G. Guseinov$^3$

\vspace{3mm}

$^1$Cumhuriyet University, Faculty of Education, Department of Mathematics and Science Education \\ 58140 Sivas, Turkey

e-mail: nhuseyin@cumhuriyet.edu.tr

\vspace{2mm}

{\small $^1$Cumhuriyet University, Faculty of Science, Department of Statistics and Computer Sciences \\ 58140 Sivas, Turkey

e-mail: ahuseyin@cumhuriyet.edu.tr

\vspace{2mm}

$^3$Eskisehir Technical  University, Faculty of Science, Department of Mathematics \\ 26470 Eskisehir, Turkey

e-mail: kguseynov@eskisehir.edu.tr}

\end{center}

\textbf{Abstract.} In this paper approximations of the set of trajectories and integral funnel of the control system described by nonlinear ordinary differential equation with integral constraint on the control functions are considered. The set of admissible control functions is replaced by a set, consisting of a finite number of piecewise-constant control functions. It is shown that the set of trajectories generated by a finite number of piecewise-constant control functions is an internal approximation of the set of trajectories. Further, each trajectory generated by a piecewise-constant control function is substituted by appropriate Euler's broken line and it is proved that the set consisting of a finite number of Euler's broken lines is an approximation of the set of trajectories of given control system. An approximation of the system's integral funnel by a set consisting of a finite number of points is obtained.

\vspace{5mm}


\vspace{5mm}


\section{Introduction}

One of the important notions of the theory of control systems described by ordinary differential equations is the attainable set and also the integral funnel concept of a given system. The attainable set is defined in the space of states  and consists of points to which the trajectories of the system arrive at the given instant of time. The integral funnel of the system  is defined in the space of positions which consists of graphs of all possible trajectories of the system, and is a generalization of the integral curve notion from the theory of ordinary differential equations. Another important notion of the theory of control systems  is the set of trajectories generated by all possible admissible control functions which includes comprehensive information about the system's behaviour. Note that if the system's trajectories are continuous functions, then the attainable set of a control system can be defined as a section of the set of trajectories at the given instant of time.

Depending on the character of the control effort, the control systems can be characterized as follows: control systems with geometric constraints on the control functions; control systems with integral constraints on the control functions; control systems with mixed constraints on the control functions, which includes both the geometric and the integral constraints. Different topological properties and methods of approximate construction of the attainable sets and the integral funnel of the control systems with geometric constraints on the control functions have been studied in numerous papers. It should also be noted the studies that were carried out in the framework of the theory of differential inclusions and differential games.

An integral constraint on the control functions is inevitable if a control resource is exhausted by consumption, say as energy, fuel, finance, food, etc.  Integral constraint on the control functions does not yield the geometric boundedness, and therefore for investigation of the control systems with integral constraint on the control functions is required to apply different approaches.  In papers \cite{kra2} -- \cite{sir} the attainable sets of the linear systems, in papers \cite{gus2} -- \cite{gusm1} approximation of  the attainable sets of the systems which are affine with respect to the control vector, and in papers \cite{gus3}, \cite{gus4} approximation of the attainable sets of the nonlinear systems are studied.  For linear and affine systems one can obtain the approximations with error evaluation. Moreover, for linear and affine systems it is possible to present approximations not only for the attainable sets, but also for the set of trajectories which is a compact subset of the space of continuous functions. It should be underlined that the method, presented in the paper \cite{gus3} for approximation of the attainable sets of the nonlinear control systems, does not allow to specify an approximation of the set of trajectories.

In this paper an approximation of the set of trajectories of the control system described by the nonlinear ordinary differential equation is considered. The control functions have an integral constraint, more precisely, a closed ball of the space $L_p$, $p>1$, centered at the origin  with radius $r$, is chosen as the set of admissible control functions. An approximation of the set of trajectories is discussed. Applying obtained result, an approximation of the integral funnel of the considered control system by a set, consisting of a finite number of points is given.

The paper is organized as follows. In Section 2, the basic conditions which satisfies the system's equation and auxiliary propositions which are used in following arguments, are given. In Section 3, step by step way, the set of admissible control functions is replaced by the set which consists of a finite number of piecewise-constant control functions and generates a finite number of trajectories. It is proved that in the appropriate settings of discretization parameters, the set consisting of a finite number of trajectories is an internal approximation of the set of trajectories of the considered control system (Theorem \ref{t1}). In Section 4, the trajectories generated by piecewise-constant control functions are substituted by appropriate Euler's broken lines and approximation of the set of trajectories by a set consisting of a finite number of the Euler's broken lines is presented (Theorem \ref{t3}). Applying this result, an approximation of the system's integral funnel by the set which consists of a finite number of points is obtained (Theorem \ref{t4}).

\section{The System's Dynamics}

Consider control system described by the  nonlinear ordinary differential equation
\begin{equation} \label{eq1}
\displaystyle \dot{x}(t)=f\left(t,x(t),u(t)\right), \ \ x(t_0)=x_0
\end{equation}
 where  $x(t)\in \mathbb{R}^n$ is the phase state vector, $u(t)\in \mathbb{R}^m$ is the control vector, $t \in [t_0,\theta]$ is the time.

For given $p>1$ and $r >0$ we denote
\begin{eqnarray*}
U_{p,r}=\left\{u(\cdot) \in L_p\big([t_0,\theta];\mathbb{R}^m\big):
\left\| u(\cdot) \right\|_p \leq r \right\}
\end{eqnarray*}
where $L_p\left([t_0,\theta];\mathbb{R}^m\right)$ is the space of Lebesgue measurable functions $u(\cdot):[t_0,\theta]\rightarrow \mathbb{R}^m$ such that
$\left\|u(\cdot)\right\|_p <+\infty,$ $\displaystyle \left\|u(\cdot)\right\|_p =\left(\int_{t_0}^{\theta} \left\| u(s)\right\|^p ds\right)^{\frac{1}{p}},$ $\left\| \cdot \right\|$ denotes the Euclidean norm.

$U_{p,r}$ is called the set of admissible control functions and every  $u(\cdot) \in U_{p,r}$ is said to be an admissible control function.

It is assumed that the right-hand side of the system (\ref{eq1}) satisfies the following conditions:

\vspace{2mm}

\textbf{2.A.} The function  $f(\cdot,\cdot,\cdot ):[t_0,\theta]\times \mathbb{R}^n\times \mathbb{R}^m\rightarrow
\mathbb{R}^{n}$ is continuous;

\vspace{2mm}

\textbf{2.B.} For every bounded set $D\subset [t_0,\theta] \times \mathbb{R}^n$ there exist $\gamma_1=\gamma_1(D)> 0$, $\gamma_2=\gamma_2(D)> 0$ and $\gamma_3=\gamma_3(D)$ such that the inequality
\begin{eqnarray*}
\left\| f(t,x_1, u_1)-f(t,x_{2},u_2)\right\| \leq \left[\gamma_1+\gamma_2 (\|u_1\|+\|u_2\|)\right] \left\| x_{1}-x_{2}\right\|+\gamma_3 \|u_1-u_2\|
\end{eqnarray*} is satisfied for every $(t,x_1,u_1)\in D\times \mathbb{R}^m$ and $(t,x_2,u_2)\in D\times \mathbb{R}^m$;

\vspace{2mm}

\textbf{2.C.} There exists $c>0$ such that the inequality
\begin{eqnarray*}
\left\| f(t,x, u)\right\| \leq c \left(\| x\| +1\right) \left(\| u \|+1\right)
\end{eqnarray*} is held for every $(t,x,u)\in [t_0,\theta] \times \mathbb{R}^n \times \mathbb{R}^m$.

\vspace{2mm}
If the function $(t,x,u)\rightarrow f(t,x,u ):[t_0,\theta]\times \mathbb{R}^n\times \mathbb{R}^m\rightarrow
\mathbb{R}^{n}$  is Lip\-schitz continuous with respect to $(x,u)$, then the conditions 2.B and 2.C are satisfied.

Let us define the trajectory of the system (\ref{eq1}) generated by an admissible control function $u_*(\cdot)\in
U_{p,r}$. An absolutely continuous function $x_*(\cdot): [t_0,\theta] \rightarrow \mathbb{R}^n$ satisfying the equation $\dot{x}_*(t)=f\left(t,x_*\left(t\right),u_*(t)\right)$ for almost all $t\in [t_0,\theta]$ and initial condition $x_*(t_0)=x_0$  is said to be a trajectory of the system (\ref{eq1}) generated by the admissible control function $u_*(\cdot)\in U_{p,r} \ .$ The set of trajectories of the system (\ref{eq1}) generated by all admissible control functions $u(\cdot)\in U_{p,r}$ is denoted by $X_{p,r}(t_0,x_0)$ and is called briefly the set of trajectories of the system (\ref{eq1}).

We set
\begin{equation}\label{ats}
X_{p,r}(t;t_0,x_0)= \left\{ x(t) \in \mathbb{R}^n: x(\cdot)\in X_{p,r}(t_0,x_0)\right\}, \ t\in [t_0,\theta],
\end{equation}
\begin{equation}\label{intfun}
F_{p,r}(t_0,x_0)= \left\{ (t,x(t))\in [t_0,\theta] \times \mathbb{R}^n: x(\cdot)\in X_{p,r}(t_0,x_0)\right\}.
\end{equation}

The set $X_{p,r}(t;t_0,x_0)$ is called the attainable set of the system (\ref{eq1}) at the instant of time $t$. It is obvious that the set $X_{p,r}(t;t_0,x_0)$ consists of points to which arrive the trajectories of the system (\ref{eq1}) at the instant of time $t$.

The set $F_{p,r}(t_0,x_0)$ is said to be the integral funnel of the system (\ref{eq1}). It consists of the graphs of all possible trajectories.

By symbol $C\left([t_0,\theta]; \mathbb{R}^n\right)$ we denote the space of continuous functions $x(\cdot):[t_0,\theta]\rightarrow \mathbb{R}^n$ with norm $\left\|x(\cdot)\right\|_C=\max \left\{\|x(t)\|: t\in [t_0,\theta]\right\}$.

$h_n(\cdot,\cdot)$ and $h_C(\cdot,\cdot)$ stand for the Hausdorff distance between the subsets of the spaces $\mathbb{R}^n$ and $C\left([t_0,\theta]; \mathbb{R}^n\right)$ respectively.

Let us formulate the propositions which will be used in following arguments.

\begin{proposition} \label{p1} \cite{gus3}  Every admissible control function $u(\cdot)\in
U_{p,r}$ generates unique trajectory of the system (\ref{eq1}).
\end{proposition}

\begin{proposition} \label{p2} \cite{gus3}
The set of trajectories $X_{p,r}(t_0,x_0)$ of the system (\ref{eq1}) is a bounded subset of the space $C\left([t_0,\theta]; \mathbb{R}^n\right)$, i.e. there exists $\alpha_*>0$ such that $\left\|x(\cdot)\right\|_C \leq \alpha_*$ for every $x(\cdot)\in X_{p,r}(t_0,x_0)$.
\end{proposition}

Let
\begin{equation}\label{fid}
\varphi (\Delta)=c(\alpha_*+1)\left( \Delta+r \Delta^{\frac{p-1}{p}}\right), \ \ \Delta >0.
\end{equation}

\begin{proposition} \label{p3} \cite{gus3}
For every $x(\cdot)\in X_{p,r}(t_0,x_0)$, $t_1\in [t_0,\theta]$ and $t_2\in [t_0,\theta]$ the inequality
\begin{eqnarray*}
\|x(t_1)-x(t_2)\| \leq \varphi \left(\left|t_1-t_2\right| \right)
\end{eqnarray*} is verified, and hence
\begin{eqnarray*}
h_n(X_{p,r}(t_1;t_0,x_0), X_{p,r}(t_2;t_0,x_0)) \leq \varphi \left(\left|t_1-t_2\right| \right)
\end{eqnarray*}
where $\varphi(\cdot)$ is defined by (\ref{fid}).
\end{proposition}

\begin{theorem} \label{t0} \cite{gus3} The set of trajectories $X_{p,r}(t_0,x_0)$ of the system (\ref{eq1}) is a precompact subset of the space $C\left([t_0,\theta]; \mathbb{R}^n\right)$.
\end{theorem}

Note that the set of trajectories of the system (\ref{eq1}) is not a closed set in general (see, \cite{gus5}). Denote
\begin{eqnarray*}
B_n(\alpha_*)= \left\{ x \in \mathbb{R}^n: \|x\| \leq \alpha_*\right\},
\end{eqnarray*}
\begin{equation} \label{den}
D_n(\alpha_*)= \left\{ (t,x) \in [t_0,\theta] \times \mathbb{R}^n: x \in B_n(\alpha_*)\right\}
\end{equation} where $\alpha_*$ is defined in Proposition \ref{p2}.

Here and henceforth we will have in mind the cylinder $D_n(\alpha_*)$ as the set $D$ in Condition 2.B.

\section{Approximation of the Set of Trajectories}

Let $\beta>0$ and $\sigma>0$ be given numbers, $\Gamma=\left\{t_0,t_1,\ldots, t_N=\theta\right\}$ be a uniform partition of the closed interval  $[t_0,\theta],$ $\Lambda=\left\{0=r_0,r_1,\ldots, r_q=\beta \right\}$ be a uniform partition of the closed interval $[0,\beta]$, $\Delta=t_{i+1}-t_i,$ $i=0,1,\ldots, N-1,$ $\delta=r_{j+1}-r_j,$ $j=0,1,\ldots, q-1,$  $S=\left\{x\in \mathbb{R}^m: \left\|x\right\|=1\right\},$ $S_{\sigma}=\left\{b_1,b_2, \ldots, b_a\right\}$ be a finite $\sigma$-net on the compact $S.$ Note that an algorithm for specifying a finite $\sigma$-net on $S$ is given in \cite{gus4}. Denote

\begin{eqnarray}\label{eq2}
U_{p,r}^{\beta, \Delta, \delta,\sigma}= \big\{ u(\cdot)\in L_p([t_0,\theta]; \mathbb{R}^m): u(t)=r_{j_i}b_{l_i} \ \mbox{for every} \ t \in [t_i,t_{i+1}), \nonumber \\ r_{j_i} \in \Lambda, \ b_{l_i}\in S_{\sigma},  \ i=0,1,\ldots, N-1, \ \Delta \cdot \sum_{i=0}^{N-1} r_{j_i}^p \leq r^p \big\}.
\end{eqnarray}

By symbol $X_{p,r}^{\beta,\Delta,\delta,\sigma}(t_0,x_0)$ we denote the set of trajectories of the system (\ref{eq1}) generated by all control functions $u(\cdot)\in U_{p,r}^{\beta,\Delta,\delta,\sigma}.$ It is obvious that the set $U_{p,r}^{\beta, \Delta, \delta,\sigma}$ consists of a finite number of piecewise-constant functions, the set $X_{p,r}^{\beta, \Delta, \delta,\sigma}(t_0,x_0)$ consists of a finite number of trajectories. For given $t\in [t_0,\theta]$ we set
\begin{equation}\label{atsap}
X_{p,r}^{\beta,\Delta,\delta,\sigma}(t;t_0,x_0)= \left\{ x(t) \in \mathbb{R}^n: x(\cdot)\in X_{p,r}^{\beta,\Delta,\delta,\sigma}(t_0,x_0)\right\}.
\end{equation}

\begin{theorem}\label{t1}
For every $\varepsilon >0$ there exist $\beta_*(\varepsilon)>0$, $\Delta_*(\varepsilon)>0$, $\delta_*(\varepsilon)>0$, $\sigma_*(\varepsilon,\beta_*(\varepsilon))>0$ such that for each $\Delta \in \left(0,\Delta_*(\varepsilon)\right],$ $\delta \in \left(0,\delta_*(\varepsilon)\right]$ and $\sigma \in (0,\sigma_*(\varepsilon,\beta_*(\varepsilon))]$, the inequality
\begin{equation}\label{os1}
h_{C}\left(X_{p,r}(t_0,x_0), X_{p,r}^{\beta_*(\varepsilon),\Delta, \delta,\sigma}(t_0,x_0)\right) \leq \frac{\varepsilon}{2}
\end{equation}
is satisfied where  $\Delta$ is the diameter of the uniform partition $\Gamma$ of the closed interval $[t_0,\theta]$,  $\delta$ is the diameter of the uniform partition $\Lambda$  of the closed interval  $[0,\beta_*(\varepsilon)]$.
\end{theorem}

\begin{proof}
For given $\beta>0$ we set
\begin{eqnarray*}
U_{p,r}^{\beta}=\left\{u(\cdot) \in U_{p,r}: \left\| u(t) \right\| \leq \beta \ \mbox{for almost all} \ t\in [t_0,\theta] \right\},
\end{eqnarray*} and let $X_{p,r}^{\beta}(t_0,x_0)$ be the set of trajectories of the system (\ref{eq1}) generated by all control functions $u(\cdot)\in U_{p,r}^{\beta}.$

Denote
\begin{equation}\label{eq3}
\kappa_*= 2\gamma_3r^p  \exp(c_0)
\end{equation} where
\begin{equation}\label{c0}
c_0=\gamma_1 (\theta -t_0) +2\gamma_2 rl_* \, ,
\end{equation}
\begin{equation}\label{el*}
l_*= \max \{ \theta-t_0, 1 \} .
\end{equation}

According to Proposition 2.1 from \cite{gus3} we have that for every $\beta>0$ the inequality
\begin{equation}\label{eq4}
h_{C}\left(X_{p,r}(t_0,x_0), X_{p,r}^{\beta}(t_0,x_0)\right) \leq \frac{\kappa_*}{\beta^{p-1}}
\end{equation}
is verified where $\kappa_*$ is defined by (\ref{eq3}). We set
\begin{equation}\label{eq5}
\beta_*(\varepsilon)=\left(\frac{10\kappa_*}{\varepsilon}\right)^{\frac{1}{p-1}}.
\end{equation}

From (\ref{eq4}) and (\ref{eq5}) it follows that
\begin{equation}\label{eq6}
h_{C}\left(X_{p,r}(t_0,x_0), X_{p,r}^{\beta_*(\varepsilon)}(t_0,x_0) \right) \leq \frac{\varepsilon}{10} \, .
\end{equation}

Now we narrow down the set of control functions $U_{p,r}^{\beta_*(\varepsilon)}$ and define new set of control functions which are Lipschitz continuous and satisfy mixed, i.e. the integral and geometric constraints. We set
\begin{equation}\label{eq7}
U_{p,r}^{\beta_*(\varepsilon),Lip}=\left\{u(\cdot) \in U_{p,r}^{\beta_*(\varepsilon)}:  u(\cdot): [t_0,\theta] \rightarrow \mathbb{R}^m \  \mbox{is Lipschitz continuous}  \right\},
\end{equation} and let $X_{p,r}^{\beta_*(\varepsilon),Lip}(t_0,x_0)$ be the set of trajectories of the system (\ref{eq1}) generated by all control functions $u(\cdot)\in U_{p,r}^{\beta_*(\varepsilon),Lip}.$

According to the Proposition 3.3 from \cite{gus3} we have that
\begin{equation}\label{eq8}
h_C\left( X_{p,r}^{\beta_*(\varepsilon)}(t_0,x_0), X_{p,r}^{\beta_*(\varepsilon),Lip}(t_0,x_0) \right)=0 \, .
\end{equation}

For given $\beta_*(\varepsilon)>0$ and integer $R>0$ let us denote
\begin{eqnarray}\label{eq9}
U_{p,r}^{\beta_*(\varepsilon), Lip, R} = \big\{ u(\cdot)\in U_{p,r}^{\beta_*(\varepsilon), Lip} : \ \mbox{the Lipschitz constant of} \ u(\cdot) \ \mbox{is not greater than} \ R \big\}
\end{eqnarray} where $U_{p,r}^{\beta_*(\varepsilon), Lip}$ is defined by (\ref{eq7}). Let $X_{p,r}^{\beta_*(\varepsilon),Lip,R}$ be the set of trajectories of the system (\ref{eq1}) generated by all control functions $u(\cdot)\in U_{p,r}^{\beta_*(\varepsilon),Lip,R}$. It is possible to show that $U_{p,r}^{\beta_*(\varepsilon), Lip, R}\subset C([t_0,\theta];\mathbb{R}^m)$ and $X_{p,r}^{\beta_*(\varepsilon), Lip, R} (t_0,x_0) \subset C([t_0,\theta];\mathbb{R}^n)$ are compact sets.

Since  $U_{p,r}^{\beta_*(\varepsilon), Lip} = \displaystyle \bigcup_{R=1}^{+\infty}U_{p,r}^{\beta_*(\varepsilon), Lip, R}$, then
\begin{equation}\label{eq10}
X_{p,r}^{\beta_*(\varepsilon),Lip}(t_0,x_0) = \bigcup_{R=1}^{+\infty} X_{p,r}^{\beta_*(\varepsilon), Lip,R}(t_0,x_0)
\end{equation} where $U_{p,r}^{\beta_*(\varepsilon), Lip, R}$ is defined by (\ref{eq9}).

Now (\ref{eq6}), (\ref{eq8}) and (\ref{eq10}) imply that the inequality
\begin{equation}\label{eq10a}
h_C\left( X_{p,r}(t_0,x_0), \bigcup_{R=1}^{+\infty} X_{p,r}^{\beta_*(\varepsilon),Lip,R}(t_0,x_0) \right) \leq \frac{\varepsilon}{10}
\end{equation} holds, where $\beta_*(\varepsilon)$ is defined by (\ref{eq5}).

According to the Theorem \ref{t0} the set of trajectories $X_{p,r}(t_0,x_0)$ is a precompact subset of the space $C([t_0,\theta];\mathbb{R}^n)$. Since
\begin{eqnarray*}
X_{p,r}^{\beta_*(\varepsilon),Lip,R}(t_0,x_0) \subset X_{p,r}^{\beta_*(\varepsilon),Lip,R+1}(t_0,x_0) \subset X_{p,r}(t_0,x_0)
\end{eqnarray*} for every $R=1,2,\ldots$, then from (\ref{eq10a}) it follows that there exists $R_*(\varepsilon)>0$ such that for any $R\geq R_*(\varepsilon)$ the inequality
\begin{eqnarray*}
h_C\left( X_{p,r}(t_0,x_0), X_{p,r}^{\beta_*(\varepsilon),Lip,R}(t_0,x_0) \right) \leq \frac{\varepsilon}{5}
\end{eqnarray*} is satisfied. In particular, the last inequality yields that the inequality
\begin{equation}\label{eq11}
h_C\left( X_{p,r}(t_0,x_0), X_{p,r}^{\beta_*(\varepsilon),Lip,R_*(\varepsilon)}(t_0,x_0) \right) \leq \frac{\varepsilon}{5}
\end{equation} is verified.

Now  we introduce new set of control functions which consists of piecewise-constant functions satisfying the integral and geometric constraints. For uniform partition $\Gamma=\left\{t_0,t_1,\ldots, t_N=\theta\right\}$  of the closed interval  $[t_0,\theta]$ where $\Delta=t_{i+1}-t_i$, $i=0,1,\ldots, N-1,$ we set
\begin{eqnarray}\label{eq12}
U_{p,r}^{\beta_*(\varepsilon), \Delta}= \big\{ u(\cdot)\in U_{p,r}^{\beta_*(\varepsilon)} : u(t)=u_i \ \mbox{for every} \ t \in [t_i,t_{i+1}), \ i=0,1,\ldots, N-1 \big\},
\end{eqnarray} and let $X_{p,r}^{\beta_*(\varepsilon),\Delta}(t_0,x_0)$ be the set of trajectories of the system (\ref{eq1}) generated by all control functions $u(\cdot)\in U_{p,r}^{\beta_*(\varepsilon),\Delta}.$

Denote
\begin{equation}\label{g1}
g_1= \gamma_3 (\theta-t_0)\cdot \exp \left( c_0\right)
\end{equation} where $c_0$ is defined by (\ref{c0}).

Let us prove that
\begin{equation}\label{eq13}
X_{p,r}^{\beta_*(\varepsilon),Lip, R_*(\varepsilon)}(t_0,x_0) \subset X_{p,r}^{\beta_*(\varepsilon),\Delta} +g_1 R_*(\varepsilon) \Delta \cdot \mathcal{B}_C(1)
\end{equation} is held where $R_*(\varepsilon)$ is defined in (\ref{eq11}),
\begin{equation}\label{eq14}
B_C(1)= \{ x(\cdot)\in C([t_0,\theta];\mathbb{R}^n) : \|x(\cdot)\|_C \leq 1\}.
\end{equation}

Choose an arbitrary $x(\cdot)\in X_{p,r}^{\beta_*(\varepsilon),Lip, R_*(\varepsilon)},$ generated by the control function $u(\cdot)\in U_{p,r}^{\beta_*(\varepsilon),Lip, R_*(\varepsilon)}.$ According to (\ref{eq9}) we have
\begin{eqnarray}\label{eq15}
\left\{
\begin{array}{llll}
\| u(\cdot)\|_p \leq r, \ \ \| u(t)\| \leq \beta_*(\varepsilon) \ \mbox{for every} \  t\in [t_0,\theta], \\
\| u(t_1)-u(t_2) \| \leq R_*(\varepsilon)\left| t_1-t_2 \right| \ \mbox{for every} \ t_1 \in [t_0,\theta], \ t_2 \in [t_0,\theta].
\end{array}
\right.
\end{eqnarray}

Define new control function $u_*(\cdot):[t_0,\theta] \rightarrow \mathbb{R}^m$, setting
\begin{eqnarray}\label{eq16}
u_*(t)=\left\{
\begin{array}{lll}
\displaystyle \frac{1}{\Delta} \int_{t_i}^{t_{i+1}} u(\tau) d\tau & \mbox{if} & t\in [t_i,t_{i+1}), \ i=0,1,\ldots, N-1, \\
\displaystyle u(t_{N-1}) & \mbox{if} & t=t_{N}
\end{array}
\right.
\end{eqnarray}
It is obvious that the function $u_*(\cdot):[t_0,\theta]\rightarrow \mathbb{R}^m$ is constant on the interval $[t_i,t_{i+1}),$ $i=0,1,\ldots, N-1,$ and $\left\| u_*(t)\right\| \leq \beta_*(\varepsilon)$ for every $t\in [t_0,\theta].$

From (\ref{eq16}) and H\"{o}lder's  inequality it follows that
\begin{eqnarray*}\label{eq17}
\int_{t_i}^{t_{i+1}}\left\|u_*(t)\right\|^p ds \leq \int_{t_i}^{t_{i+1}} \left\|u(\tau)\right\|^p d\tau
\end{eqnarray*} for every $i=0,1,\ldots, N-1,$ and hence $\left\|u_*(\cdot)\right\|_p \leq \left\|u(\cdot)\right\|_p \leq r.$
Thus we obtain that $u_*(\cdot)\in U_{p,r}^{\beta_*(\varepsilon),\Delta}$ where $U_{p,r}^{\beta_*(\varepsilon), \Delta}$ is defined by (\ref{eq12}). Let $x_*(\cdot)$ be the trajectory of the system (\ref{eq1}) generated by the control function $u_*(\cdot).$ Then $x_*(\cdot)\in X_{p,r}^{\beta_* (\varepsilon), \Delta}(t_0,x_0).$

Choose an arbitrary $t_*\in [t_0,\theta]$ and fix it. Then there exists $i_*=0,1,\ldots, N-1$ such that $t_*\in [t_{i_*},t_{i_*+1}).$ If $t_*=\theta$, then we have $t_*\in [t_{N-1},\theta].$ From (\ref{eq15}), (\ref{eq16}) and equality $t_{i_*+1}-t_{i_*} =\Delta$ it follows that
\begin{eqnarray*}
\left\| u(t_*)-u_*(t_*)\right\| \leq \frac{1}{\Delta} \int_{t_{i_*}}^{t_{i_*+1}}\left\| u(t_*)- u(\tau)\right\| d\tau \leq
\frac{1}{\Delta} \int_{t_{i_*}}^{t_{i_*+1}}R_*(\varepsilon)\left| t_*- \tau \right| d\tau \leq R_*(\varepsilon)\Delta.
\end{eqnarray*}
Since $t_*\in [t_0,\theta]$ is arbitrarily fixed, then we have that
\begin{equation}\label{eq18}
\left\| u(t)-u_*(t)\right\| \leq  R_*(\varepsilon)\Delta
\end{equation} for every $t\in [t_0,\theta]$.

Now, from Conditions  2.B and (\ref{eq18}) it follows that
\begin{eqnarray}\label{eq19}
\left\| x(t)-x_*(t) \right\| \leq \int_{t_0}^{t} \big[\gamma_1 +  \gamma_2(\|u(\tau)\|+\|u_*(\tau)\|\big]\cdot \| x(\tau)-x_*(\tau) \| d\tau + \gamma_3 (\theta-t_0)R_*(\varepsilon)\Delta
\end{eqnarray} for every $t\in [t_0,\theta]$. Since $u(\cdot)\in U_{p,r}^{\beta_*(\varepsilon),Lip, R_*(\varepsilon)}\subset U_{p,r},$ $u_*(\cdot)\in U_{p,r}^{\beta_*(\varepsilon),\Delta} \subset U_{p,r},$ then the Gronwall-Bellman inequality, (\ref{c0}), (\ref{el*}), (\ref{g1}) and (\ref{eq19}) imply that
\begin{eqnarray}\label{eq20}
\left\| x(t)-x_*(t) \right\| & \leq & \gamma_3 (\theta-t_0)R_*(\varepsilon)\Delta \cdot  \exp \left(\int_{t_0}^{\theta} \big[\gamma_1 +  \gamma_2(\|u(\tau)\|+\|u_*(\tau)\|\big] d \tau \right) \nonumber \\ & \leq & \gamma_3 (\theta-t_0)R_*(\varepsilon)\Delta \cdot \exp (c_0)  =g_1 R_*(\varepsilon) \Delta
\end{eqnarray} for every  $t\in [t_0,\theta]$. Since  $x(\cdot)\in X_{p,r}^{\beta_*(\varepsilon),Lip, R_*(\varepsilon)}(t_0,x_0)$ is arbitrarily chosen,  $x_*(\cdot)\in X_{p,r}^{\beta_*(\varepsilon),\Delta}(t_0,x_0)$, then form (\ref{eq20}) we obtain the validity of the inclusion (\ref{eq13}).

Denote
\begin{equation}\label{eq21}
\Delta_*(\varepsilon)= \frac{\varepsilon}{10g_1 R_*(\varepsilon)} \, .
\end{equation}

The inclusion (\ref{eq13}) and the equality (\ref{eq21}) imply that for every uniform partition $\Gamma= $ \linebreak $\{t_0,t_1,\ldots, t_N=\theta \}$  of the closed interval $[t_0,\theta]$ such that $\Delta \leq \Delta_*(\varepsilon),$ the inclusion
\begin{equation}\label{eq22}
X_{p,r}^{\beta_*(\varepsilon),Lip, R_*(\varepsilon)}(t_0,x_0) \subset X_{p,r}^{\beta_*(\varepsilon),\Delta}(t_0,x_0) + \frac{\varepsilon}{10} \mathcal{B}_C(1)
\end{equation} is satisfied where $\mathcal{B}_C(1)$ is defined by (\ref{eq14}), $\Delta=t_{i+1}-t_i$, $i=0,1,\ldots, N-1$.

On behalf of the inequality (\ref{eq11}) and inclusion (\ref{eq22}) we have that for every uniform partition $\Gamma=\left\{t_0,t_1,\ldots, t_N=\theta\right\}$  of the closed interval $[t_0,\theta]$ such that $\Delta \leq \Delta_*(\varepsilon),$ the inclusion
\begin{equation}\label{eq23}
X_{p,r}(t_0,x_0) \subset X_{p,r}^{\beta_*(\varepsilon),\Delta}(t_0,x_0) + \frac{3\varepsilon}{10} \cdot \mathcal{B}_C(1)
\end{equation} is verified.

Since $X_{p,r}^{\beta_*(\varepsilon),\Delta}(t_0,x_0) \subset X_{p,r}(t_0,x_0)$, then from (\ref{eq23}) we obtain that for every uniform partition $\Gamma=\left\{t_0,t_1,\ldots, t_N=\theta\right\}$  of the closed interval $[t_0,\theta]$ such that $\Delta \leq \Delta_*(\varepsilon),$ the inequality
\begin{equation}\label{eq24}
h_C \left(X_{p,r}(t_0,x_0), X_{p,r}^{\beta_*(\varepsilon),\Delta}(t_0,x_0)\right) \leq \frac{3\varepsilon}{10}
\end{equation} is held.

For given $\beta_*(\varepsilon),$ uniform partition $\Gamma=\left\{t_0,t_1,\ldots, t_N=\theta\right\}$ of the closed interval $[t_0,\theta]$ and uniform partition $\Lambda=\left\{0=r_0,r_1,\ldots, r_q=\beta_*(\varepsilon)\right\}$ of the closed interval $[0,\beta_*(\varepsilon)]$ where $\Delta =t_{i+1}-t_i,$ $i=0,1,\ldots, N-1,$ is diameter of the uniform partition $\Gamma,$   $\delta =r_{j+1}-r_j,$ $j=0,1,\ldots, q-1,$ is diameter of the uniform partition $\Lambda,$ we define new set of control functions, setting
\begin{eqnarray}\label{eq25}
U_{p,r}^{\beta_*(\varepsilon), \Delta, \delta}&=& \big\{ u(\cdot)\in U_{p,r}^{\beta_*(\varepsilon),\Delta} : u(s)=u_i \ \mbox{for every} \ t \in [t_i,t_{i+1}) \nonumber \\ & & \quad \mbox{and} \ \left\|u_i\right\| \in \Lambda,  \ i=0,1,\ldots, N-1  \big\}.
\end{eqnarray} Let $X_{p,r}^{\beta_*(\varepsilon),\Delta,\delta}(t_0,x_0)$ be  the set of trajectories of the system (\ref{eq1}) generated by all control functions $u(\cdot)\in U_{p,r}^{\beta_*(\varepsilon),\Delta,\delta}.$

According to the Proposition 6.1 from \cite{gus3} we have that for every uniform partition $\Gamma=$ \linebreak $\left\{t_0,t_1,\ldots, t_N=\theta\right\}$ of the interval $[t_0,\theta]$ and uniform partition $\Lambda=\left\{0=r_0,r_1,\ldots, r_q=\beta_*(\varepsilon)\right\}$ of the closed interval $[0,\beta_*(\varepsilon)],$ the inequality
\begin{equation}\label{eq26}
h_C \left(X_{p,r}^{\beta_*(\varepsilon),\Delta}(t_0,x_0), X_{p,r}^{\beta_*(\varepsilon),\Delta,\delta}(t_0,x_0)\right) \leq g_1 \delta
\end{equation} is held where $g_1$ is defined by (\ref{g1}), $\Delta=t_{i+1}-t_i$, $i=0,1,\ldots, N-1$, $\delta=r_{j+1}-r_j$, $j=0,1,\ldots, q-1$.

Denote
\begin{equation}\label{eq34}
\delta_*(\varepsilon)=  \frac{\varepsilon}{10g_1} \, .
\end{equation}

By virtue of the inequalities (\ref{eq24}), (\ref{eq26}) and  (\ref{eq34}) we have that for every uniform partition $\Gamma=\left\{t_0,t_1,\ldots, t_N=\theta\right\}$  of the interval $[t_0,\theta]$ and for every uniform partition $\Lambda=\left\{0=r_0,r_1,\ldots, r_q=\beta_*(\varepsilon)\right\}$  of the closed interval $[0,\beta_*(\varepsilon)]$ such that $\Delta \leq \Delta_*(\varepsilon),$ $\delta \leq \delta_*(\varepsilon),$ the inequality
\begin{equation}\label{eq36}
h_C \left(X_{p,r}(t_0,x_0),X_{p,r}^{\beta_*(\varepsilon),\Delta,\delta}(t_0,x_0)\right) \leq \frac{2\varepsilon}{5}
\end{equation} is verified where $\Delta=t_{i+1}-t_i$, $i=0,1,\ldots, N-1$, $\delta=r_{j+1}-r_j$, $j=0,1,\ldots, q-1$.

On behalf of the Proposition 7.1 from \cite{gus3} we obtain that for every uniform partition $\Gamma= $ \linebreak $\left\{t_0,t_1,\ldots, t_N=\theta\right\}$ of the interval $[t_0,\theta]$ and uniform partition $\Lambda=\left\{0=r_0,r_1,\ldots, r_q=\beta_*(\varepsilon)\right\}$ of the closed interval $[0,\beta_*(\varepsilon)]$ and $\sigma >0$ the inequality
\begin{equation}\label{eq37}
h_C \left(X_{p,r}^{\beta_*(\varepsilon),\Delta,\delta}(t_0,x_0), X_{p,r}^{\beta_*(\varepsilon),\Delta,\delta,\sigma}(t_0,x_0)\right) \leq g_1 \beta_*(\varepsilon)\sigma
\end{equation} holds where $\Delta =t_{i+1}-t_i,$ $j=0,1,\ldots, N-1,$ is diameter of the uniform partition $\Gamma,$ $\delta =r_{j+1}-r_j,$ $j=0,1,\ldots, q-1,$ is diameter of the uniform partition $\Lambda,$ $g_1$ is defined by (\ref{g1}), $X_{p,r}^{\beta_*(\varepsilon),\Delta,\delta}(t_0,x_0)$ and $X_{p,r}^{\beta_*(\varepsilon),\Delta,\delta,\sigma}(t_0,x_0)$
are the sets of trajectories of the system (\ref{eq1}), generated by the set of control functions $U_{p,r}^{\beta_*(\varepsilon),\Delta,\delta}$ and $U_{p,r}^{\beta_*(\varepsilon),\Delta,\delta,\sigma}$ respectively. The set $U_{p,r}^{\beta_*(\varepsilon),\Delta,\delta}$ is defined by (\ref{eq25}), $U_{p,r}^{\beta_*(\varepsilon),\Delta,\delta,\sigma}$ is defined by (\ref{eq2}).

Let us set
\begin{equation}\label{eq45}
\sigma_*(\varepsilon,\beta_*(\varepsilon)) =\frac{\varepsilon}{10 g_1 \beta_*(\varepsilon)} \, .
\end{equation}

By virtue of (\ref{eq36}), (\ref{eq37}) and  (\ref{eq45}) we have that for every uniform partition $\Gamma=\left\{t_0,t_1,\ldots, t_N=\theta\right\}$  of the interval $[t_0,\theta]$, for every uniform partition $\Lambda=\left\{0=r_0,r_1,\ldots, r_q=\beta_*(\varepsilon)\right\}$  of the closed interval $[0,\beta_*(\varepsilon)]$ and for every $\sigma >0$ such that $\Delta \leq \Delta_*(\varepsilon),$ $\delta \leq \delta_*(\varepsilon),$ $\sigma \leq \sigma_*(\varepsilon, \beta_*(\varepsilon))$, the inequality (\ref{os1}) is verified where $\beta_*(\varepsilon),$ $\Delta_*(\varepsilon),$ $\delta_*(\varepsilon)$ and $\sigma_*(\varepsilon,\beta_*(\varepsilon))$ are defined by (\ref{eq5}), (\ref{eq21}), (\ref{eq34}) and (\ref{eq45}) respectively, $\Delta=t_{i+1}-t_i$, $i=0,1,\ldots, N-1$, $\delta=r_{j+1}-r_j$, $j=0,1,\ldots, q-1$.

The theorem  is proved.
\end{proof}

From Theorem \ref{t1} we obtain the validity of the following corollary which gives us an approximation of the attainable sets.

\begin{corollary}\label{cor1}
For every $\varepsilon >0$ and  for each $\Delta \in \left(0,\Delta_*(\varepsilon)\right],$ $\delta \in \left(0,\delta_*(\varepsilon)\right]$ and $\sigma \in (0,\sigma_*(\varepsilon,\beta_*(\varepsilon))]$, the inequality
\begin{eqnarray*}
h_{n}\left(X_{p,r}(t;t_0,x_0), X_{p,r}^{\beta_*(\varepsilon),\Delta, \delta,\sigma}(t;t_0,x_0)\right) \leq \frac{\varepsilon}{2}
\end{eqnarray*}
is verified where  $\Delta$ is the diameter of the uniform partition $\Gamma$ of the closed interval $[t_0,\theta]$,  $\delta$ is the diameter of the uniform partition $\Lambda$  of the closed interval  $[0,\beta_*(\varepsilon)]$, $\beta_*(\varepsilon),$ $\Delta_*(\varepsilon),$ $\delta_*(\varepsilon)$ and $\sigma(\varepsilon,\beta_*(\varepsilon))$ are defined by (\ref{eq5}), (\ref{eq21}), (\ref{eq34}) and (\ref{eq45}) respectively.

Here $X_{p,r}(t;t_0,x_0)$ is attainable set of the system (\ref{eq1}) at the instant of time $t$ and is defined by (\ref{ats}),
$X_{p,r}^{\beta_*(\varepsilon),\Delta, \delta,\sigma}(t;t_0,x_0)$ is defined by (\ref{atsap}) and consists of a finite number of points.
\end{corollary}

\begin{remark} \label{r1}
Since $X_{p,r}^{\beta_*(\varepsilon),\Delta,\delta,\sigma}(t_0,x_0) \subset X_{p,r}(t_0,x_0)$, then we obtain that the presented approximation is an internal one.
\end{remark}

\section{Euler's Broken Lines and Approximation of the Integral Funnel}

For given $\beta_*(\varepsilon)>0$,  uniform partition $\Gamma=\left\{t_0,t_1,\ldots, t_N=\theta\right\}$ of the closed interval  $[t_0,\theta],$ uniform partition $\Lambda=\left\{0=r_0,r_1,\ldots, r_q=\beta \right\}$ of the closed interval $[0,\beta]$ and a finite $\sigma$-net $S_{\sigma}=\left\{b_1,b_2, \ldots, b_a\right\}$ we set
\begin{eqnarray}\label{eq48}
&& Z_{p,r}^{\beta_*(\varepsilon), \Delta, \delta,\sigma}(t_0,x_0) =  \big\{ z(\cdot)\in C([t_0,\theta]; \mathbb{R}^n): z(t)=z(t_i) +(t-t_i)f(t_i,x(t_i),r_{j_i}b_{l_i}), \nonumber \\ && \quad t \in [t_i,t_{i+1}], \  z(t_0)=x_0,  \  r_{j_i} \in \Lambda, \ b_{l_i}\in S_{\sigma},  \ i=0,1,\ldots, N-1, \ \Delta \cdot \sum_{i=0}^{N-1} r_{j_i}^p \leq r^p \big\},
\end{eqnarray}
\begin{equation}\label{eq48*}
Z_{p,r}^{\beta_*(\varepsilon), \Delta, \delta,\sigma}(t;t_0,x_0)= \big\{ z(t)\in \mathbb{R}^n: z(\cdot)\in Z^{\beta_*(\varepsilon), \Delta, \delta,\sigma}(t_0,x_0) \big\}
\end{equation}
where $\beta_*(\varepsilon)>0$ is defined by (\ref{eq5}), $\Delta=t_{i+1}-t_i,$ $i=0,1,\ldots, N-1,$ $\delta=r_{j+1}-r_j,$ $j=0,1,\ldots, q-1,$ $t\in [t_0,\theta],$  Now we denote
\begin{equation}\label{eq49}
\varphi_*(\Delta)= \max \left\{ \Delta, \varphi(\Delta))\right\},
\end{equation}
\begin{eqnarray*}
V(\alpha_*,\beta_*(\varepsilon)) = D_n(\alpha_*) \times B_m(\beta_*(\varepsilon)),
\end{eqnarray*}
\begin{eqnarray}\label{eq51}
&& \omega \left(\varphi_* (\Delta),\beta_*(\varepsilon)\right) =\max \Big\{ \|f(t_1,x_1,u)-f(t_2,x_2,u)\|:  (t_1,x_1,u)\in V(\alpha_*,\beta_*(\varepsilon)), \nonumber \\ && \qquad (t_2,x_2,u)\in V(\alpha_*,\beta_*(\varepsilon)), \
|t_1-t_2| \leq \varphi_*(\Delta), \ \|x_1-x_2\| \leq \varphi_*(\Delta) \Big\},
\end{eqnarray}
\begin{equation}\label{geb}
g(\beta_*(\varepsilon))=\gamma_1+2\gamma_2 \beta_*(\varepsilon)
\end{equation}
where $\varphi (\cdot)$ is defined by (\ref{fid}), $\beta_*(\varepsilon)$ is defined by (\ref{eq5}), $D_n(\alpha_*)$ is defined by (\ref{den}), $ B_m(\beta_*(\varepsilon))=\{ u\in \mathbb{R}^m: \| u\| \leq \beta_*(\varepsilon)\}$. For every fixed $\varepsilon >0$ we have that $\varphi_*(\Delta)\rightarrow 0$ and $\omega \left(\varphi_* (\Delta),\beta_*(\varepsilon)\right) \rightarrow 0$ as $\Delta \rightarrow 0^+$.

\begin{proposition}\label{p4.1}
For every $\varepsilon >0$, uniform partition $\Gamma=\left\{t_0,t_1,\ldots, t_N=\theta\right\}$ of the interval $[t_0,\theta]$,  uniform partition $\Lambda=\left\{0=r_0,r_1,\ldots, r_q=\beta_*(\varepsilon)\right\}$ of the closed interval $[0,\beta_*(\varepsilon)]$ and $\sigma >0$ the inequality
\begin{eqnarray*}
h_C \left(X_{p,r}^{\beta_*(\varepsilon),\Delta,\delta,\sigma}(t_0,x_0), Z_{p,r}^{\beta_*(\varepsilon),\Delta,\delta,\sigma}(t_0,x_0)\right) \leq \omega(\varphi_*(\Delta), \beta_*(\varepsilon))(\theta-t_0) \exp [g( \beta_*(\varepsilon))(\theta-t_0)]
\end{eqnarray*} is satisfied where $\Delta =t_{i+1}-t_i,$ $i=0,1,\ldots, N-1,$ is diameter of the uniform partition $\Gamma,$ $\delta =r_{j+1}-r_j,$ $j=0,1,\ldots, q-1,$ is diameter of the uniform partition $\Lambda,$ $\beta_*(\varepsilon)$ is defined by (\ref{eq5}), $X_{p,r}^{\beta_*(\varepsilon),\Delta,\delta,\sigma}(t_0,x_0)$ is the set of trajectories of the system (\ref{eq1}), generated by the set of control functions $U_{p,r}^{\beta_*(\varepsilon),\Delta,\delta,\sigma}$, the set $Z_{p,r}^{\beta_*(\varepsilon),\Delta,\delta,\sigma}$ is defined by (\ref{eq48}), $\omega(\varphi_*(\Delta), \beta_*(\varepsilon))$ is defined by (\ref{eq51}) and $g( \beta_*(\varepsilon))$ is defined by (\ref{geb}).
\end{proposition}

Since $\varphi_*(\Delta)\rightarrow 0$, $\omega (\varphi_*(\Delta), \beta_*(\varepsilon))\rightarrow 0$ as $\Delta \rightarrow 0^+$, then for given $\frac{\varepsilon}{10}$ there exists $\Delta^*(\varepsilon,\beta(\varepsilon))>0$ such that for every $\Delta \leq \Delta^*(\varepsilon,\beta(\varepsilon))$ the inequalities
\begin{equation}\label{eq68}
\omega (\varphi_*(\Delta), \beta_*(\varepsilon)) \leq \frac{\varepsilon}{10 \exp (g( \beta_*(\varepsilon) (\theta-t_0))}, \ \ \varphi_*(\Delta) \leq \frac{\varepsilon}{10}
\end{equation} are satisfied.

Denote
\begin{equation}\label{eq69}
\Delta_0(\varepsilon,\beta_*(\varepsilon))= \min \left\{ \Delta_*(\varepsilon), \, \Delta^*(\varepsilon, \beta_*(\varepsilon)), \,  \frac{\varepsilon}{10}\right\}
\end{equation} where $\Delta_*(\varepsilon)$ is defined by (\ref{eq21}).

Theorem \ref{t1}, Proposition \ref{p4.1},  (\ref{eq68}) and  (\ref{eq69}) yield the validity of the following theorem which presents an approximation of the set of trajectories of the system (\ref{eq1}) by a finite number of appropriate Euler's broken lines.
\begin{theorem}\label{t3}
For every $\varepsilon >0$ and for every uniform partition $\Gamma=\left\{t_0,t_1,\ldots, t_N=\theta\right\}$  of the interval $[t_0,\theta]$, for every uniform partition $\Lambda=\left\{0=r_0,r_1,\ldots, r_q=\beta_*(\varepsilon)\right\}$  of the closed interval $[0,\beta_*(\varepsilon)]$ and for every $\sigma >0$ such that $\Delta \leq \Delta_0(\varepsilon,  \beta_*(\varepsilon)),$ $\delta \leq \delta_*(\varepsilon),$ $\sigma \leq \sigma_*(\varepsilon, \beta_*(\varepsilon))$, the inequality
\begin{eqnarray*}
h_C \left(X_{p,r}(t_0,x_0),Z_{p,r}^{\beta_*(\varepsilon),\Delta,\delta,\sigma}(t_0,x_0)\right) \leq \frac{3\varepsilon}{5}
\end{eqnarray*} is verified where $\beta_*(\varepsilon),$ $\Delta_0(\varepsilon, \beta_*(\varepsilon)),$ $\delta_*(\varepsilon)$ and $\sigma_*(\varepsilon,\beta_*(\varepsilon))$ are defined by (\ref{eq5}), (\ref{eq69}), (\ref{eq34}) and (\ref{eq45}) respectively,
$\Delta =t_{i+1}-t_i,$ $i=0,1,\ldots, N-1,$  $\delta =r_{j+1}-r_j,$ $j=0,1,\ldots, q-1.$
\end{theorem}

Note that the Theorem \ref{t3} is a generalization of the main result of the paper \cite{gus3}, Theorem 9.1, where an approximation of the attainable set of the system (\ref{eq1}) is obtained.

Denote
\begin{equation}\label{eq71}
\Phi_{p,r,\Gamma}^{\beta_*(\varepsilon),\Delta,\delta,\sigma}(t_0,x_0)= \bigcup_{i=0}^{N} \left( t_i, Z_{p,r}^{\beta_*(\varepsilon),\Delta,\delta,\sigma}(t_i;t_0,x_0)\right\}
\end{equation} where the set $Z_{p,r}^{\beta_*(\varepsilon),\Delta,\delta,\sigma}(t_i;t_0,x_0)$ is defined by (\ref{eq48*}).

\begin{theorem}\label{t4}
For every $\varepsilon >0$ and for every uniform partition $\Gamma=\left\{t_0,t_1,\ldots, t_N=\theta\right\}$  of the interval $[t_0,\theta]$, for every uniform partition $\Lambda=\left\{0=r_0,r_1,\ldots, r_q=\beta_*(\varepsilon)\right\}$  of the closed interval $[0,\beta_*(\varepsilon)]$ and for every $\sigma >0$ such that $\Delta \leq \Delta_0(\varepsilon,\beta_*(\varepsilon)),$ $\delta \leq \delta_*(\varepsilon),$ $\sigma \leq \sigma_*(\varepsilon, \beta_*(\varepsilon))$, the inequality
\begin{eqnarray*}
h_{n+1} \left(F_{p,r}(t_0,x_0),\Phi_{p,r,\Gamma}^{\beta_*(\varepsilon),\Delta,\delta,\sigma}(t_0,x_0)\right) \leq \varepsilon
\end{eqnarray*} is held where $\beta_*(\varepsilon),$ $\Delta_0(\varepsilon,\beta_*(\varepsilon)),$ $\delta_*(\varepsilon)$ and $\sigma_*(\varepsilon,\beta_*(\varepsilon))$ are defined by (\ref{eq5}), (\ref{eq69}), (\ref{eq34}) and (\ref{eq45}) respectively, $\Delta =t_{i+1}-t_i,$ $i=0,1,\ldots, N-1,$  $\delta =r_{j+1}-r_j,$ $j=0,1,\ldots, q-1$, the set $F_{p,r}(t_0,x_0)$ is defined by (\ref{intfun}), the set $\Phi_{p,r,\Gamma}^{\beta_*(\varepsilon),\Delta,\delta,\sigma}(t_0,x_0)$ is defined by (\ref{eq71}).
\end{theorem}

\begin{proof}
Let $(t_*,x_*) \in F_{p,r}(t_0,x_0).$ Then $x_*\in X_{p,r}(t_*;t_0,x_0)$, and therefore there exists $x_*(\cdot)\in X_{p,r}(t_0,x_0)$ such that $x_*(t_*)=x_*.$  By virtue of Theorem \ref{t3} there exists $z_*(\cdot)\in Z_{p,r}^{\beta_*(\varepsilon),\Delta,\delta,\sigma}(t_0,x_0)$ such that
\begin{equation}\label{eq731}
\|x_*(t)-z_*(t) \|  \leq \frac{7\varepsilon}{10}
\end{equation} for every $t\in [t_0,\theta]$.

Let $t_* \in \Gamma$, i.e. let there exists $t_{i_*}\in \Gamma$ such that $t_*=t_{i_*}$. Since $z_*(\cdot)\in Z_{p,r}^{\beta_*(\varepsilon),\Delta,\delta,\sigma}(t_0,x_0)$, then (\ref{eq71}) implies that
\begin{equation}\label{eq732}
(t_*,z_*(t_*)) = (t_{i_*},z_*(t_{i_*})) \in  \Phi_{p,r,\Gamma}^{\beta_*(\varepsilon),\Delta,\delta,\sigma}(t_0,x_0).
\end{equation}

From (\ref{eq731}) it follows that
\begin{eqnarray*}
\|x_* -z_*(t_*)\| = \|x_*(t_*)-z_*(t_*)\|= \|x_*(t_{i_*})-z_*(t_{i_*}) \|  \leq \frac{7\varepsilon}{10}
\end{eqnarray*} and hence
\begin{equation}\label{eq733}
\|(t_*,x_*) -(t_*,z_*(t_*))\| = \|x_*-z_*(t_*)\|  \leq \frac{7\varepsilon}{10}
\end{equation}

Since $(t_*,x_*) \in F_{p,r}(t_0,x_0)$ is arbitrarily chosen, then (\ref{eq732}) and (\ref{eq733}) yield that
\begin{equation}\label{eq734}
F_{p,r}(t_0,x_0) \subset \Phi_{p,r,\Gamma}^{\beta_*(\varepsilon),\Delta,\delta,\sigma}(t_0,x_0)+  \frac{7\varepsilon}{10} \cdot B_{n+1}(1) \, .
\end{equation}

Now suppose that $t_* \not \in \Gamma$. Then there exists $i_0=0,1,\ldots, N-1$ such that $t_*\in (t_{i_0},t_{i_0+1}).$  Since $|t_*-t_{i_0}| \leq \Delta \leq \Delta_0(\varepsilon, \beta_*(\varepsilon)),$ then according to (\ref{eq49}), (\ref{eq68}), (\ref{eq69}) and Proposition \ref{p3} we have that
\begin{equation}\label{eq73}
\|x_*-x_*(t_{i_0}) \| =\|x_*(t_*)-x_*(t_{i_0}) \|  \leq \varphi(\Delta) \leq \varphi_*(\Delta) \leq \frac{\varepsilon}{10} \, .
\end{equation}

By virtue of (\ref{eq731}) we have that
\begin{equation}\label{eq74}
\|x_*(t_{i_0})-z_*(t_{i_0}) \|  \leq \frac{7\varepsilon}{10}
\end{equation} where $z_*(t_{i_0}) \in Z_{p,r}^{\beta_*(\varepsilon),\Delta,\delta,\sigma}(t_{i_0};t_0,x_0)$.
It is obvious that $(t_{i_0},z_*(t_{i_0})) \in \Phi_{p,r,\Gamma}^{\beta_*(\varepsilon),\Delta,\delta,\sigma}(t_0,x_0)$.
From (\ref{eq69}), (\ref{eq73}) and (\ref{eq74}) it follows that
\begin{eqnarray}\label{eq75}
\|(t_*,x_*)-(t_{i_0},z_*(t_{i_0})) \|  & \leq & |t_*-t_{i_0}| +\|x_*(t_*)-x_*(t_{i_0}) \|  +  \|x_*(t_{i_0})-z_*(t_{i_0}) \| \nonumber \\  & \leq & \frac{\varepsilon}{10}+ \frac{\varepsilon}{10}+\frac{7\varepsilon}{10} =\frac{9\varepsilon}{10} <\varepsilon \, .
\end{eqnarray}

Since $(t_*,x_*) \in F_{p,r}(t_0,x_0)$ is arbitrarily chosen, $(t_{i_0},z_*(t_{i_0})) \in \Phi_{p,r,\Gamma}^{\beta_*(\varepsilon),\Delta,\delta,\sigma}(t_0,x_0)$, then
from (\ref{eq75}) we obtain
\begin{equation}\label{eq76}
F_{p,r}(t_0,x_0) \subset \Phi_{p,r,\Gamma}^{\beta_*(\varepsilon),\Delta,\delta,\sigma}(t_0,x_0)+ \frac{9\varepsilon}{10} \cdot B_{n+1}(1) \, .
\end{equation}

Finally, (\ref{eq734}) and (\ref{eq76}) imply that
\begin{equation}\label{eq77}
F_{p,r}(t_0,x_0) \subset \Phi_{p,r,\Gamma}^{\beta_*(\varepsilon),\Delta,\delta,\sigma}(t_0,x_0)+ \varepsilon \cdot B_{n+1}(1) \, .
\end{equation}

Now let $(t_{i_*},z_0)\in  \Phi_{p,r,\Gamma}^{\beta_*(\varepsilon),\Delta,\delta,\sigma}(t_0,x_0).$ Then $z_0 \in Z_{p,r}^{\beta_*(\varepsilon),\Delta,\delta,\sigma}(t_{i_*};t_0,x_0)$ and there exists $z_0(\cdot)\in  Z_{p,r}^{\beta_*(\varepsilon),\Delta,\delta,\sigma}(t_0,x_0)$ such that $z_0(t_{i_*})=z_0$. According to the Theorem \ref{t3} there exists $x_0(\cdot)\in  X_{p,r}^{\beta_*(\varepsilon),\Delta,\delta,\sigma}(t_0,x_0)$ such that
\begin{eqnarray*}
\|x_0(t)-z_0(t) \|  \leq \frac{7\varepsilon}{10}
\end{eqnarray*} for every $t\in [t_0,\theta]$ and hence
\begin{equation}\label{eq78}
\|x_0(t_{i_*})-z_0(t_{i_*}) \|  \leq \frac{7\varepsilon}{10} \, .
\end{equation}

Since $x_0(\cdot)\in  X_{p,r}^{\beta_*(\varepsilon),\Delta,\delta,\sigma}(t_0,x_0)$ and $X_{p,r}^{\beta_*(\varepsilon),\Delta,\delta,\sigma}(t_0,x_0)\subset X_{p,r}(t_0,x_0)$, then we have $(t_{i_*},x_0(t_{i_*})) \in F_{p,r}(t_0,x_0)$. Now, from (\ref{eq78}) we obtain
\begin{equation}\label{eq79}
\|(t_{i_*},x_0(t_{i_*}))-(t_{i_*},z_0(t_{i_*})) \| = \|x_0(t_{i_*})-z_0(t_{i_*}) \| \leq \frac{7\varepsilon}{10} \, .
\end{equation}

Since $(t_{i_*},z_0)\in  \Phi_{p,r,\Gamma}^{\beta_*(\varepsilon),\Delta,\delta,\sigma}(t_0,x_0)$ is arbitrarily chosen,  $(t_{i_*},x_0(t_{i_*})) \in F_{p,r}(t_0,x_0)$, then  (\ref{eq79}) imples that
\begin{equation}\label{eq80}
\Phi_{p,r,\Gamma}^{\beta_*(\varepsilon),\Delta,\delta,\sigma}(t_0,x_0)\subset F_{p,r}(t_0,x_0)+ \varepsilon \cdot B_{n+1}(1) \, .
\end{equation}

(\ref{eq77}) and (\ref{eq80}) complete the proof of the theorem.
\end{proof}

\section*{Conclusion}

The paper presents an approximation of the set of trajectories of the nonlinear control system with integral constraint on the control functions by a finite number of Euler's broken lines, each of them is generated by a finite number of piecewise-constant control functions. Using an algorithm for specifying and aligning of the elements of a finite $\sigma$-net on the sphere $S=\left\{x\in \mathbb{R}^m: \left\|x\right\|=1\right\}$, as well as the aligning a finite number of functions from the set of control functions $U_{p,r}^{\beta_*(\varepsilon), \Delta, \delta,\sigma}$ presented in the paper \cite{gus4}, it is possible to carry out an approximate construction of the set of trajectories, attainable sets and integral funnel of the given nonlinear control system.

\end{document}